\documentclass[12pt]{article}

\usepackage{graphicx}
\usepackage{psfrag}
\usepackage{epsf}
\usepackage{amsmath,amsfonts,amssymb,latexsym,amsthm}
\usepackage{algorithmic}
\usepackage{algorithm}
\usepackage[width=160mm,height=240mm,left=35mm,foot=10mm]{geometry}
\usepackage{setspace}
\usepackage{enumerate}
\usepackage{tikz}

\newcommand{\ignore}[1]{}
\newcommand{\startClaims}{\setcounter{claim}{0}}
\newtheorem{theorem}{Theorem}[section]
\newtheorem{corollary}[theorem]{Corollary}
\newtheorem{lemma}[theorem]{Lemma}

\newtheorem{conjecture}[theorem]{Conjecture}

\title{Efficient $(j,k)$-Domination in Regular Graphs}

\author{Brendan Rooney\\
\small{Rochester Institute of Technology}}

\date{\today}

\newcommand{\eopf}{\raisebox{0.8ex}{\framebox{}}}

{\noindent{\bf Sketch of a proof.}\ }%
{\hfill\eopf\par\bigskip}%
{\noindent{\bf Ideas for a proof.}\ }%
{\hfill\eopf\par\bigskip}%
{\noindent{\bf Proof in progress.}\ }%
{\hfill\eopf\par\bigskip}%

\newcommand{\ZZ}{\mathbb Z}

\newcommand{\RR}{\mathbb R}

\newcommand{\col}{\operatorname{col}}
\newcommand{\Span}{\operatorname{span}}

\begin{document}

\maketitle

\begin{abstract}
Rubalcaba and Slater (Robert R. Rubalcaba and Peter J. Slater. Efficient (j, k)-domination. \emph{Discuss. Math. Graph Theory}, 27(3):409–423, 2007.) define a \emph{$(j,k)$-dominating function} on graph $X$ as a function $f:V(X)\rightarrow \{0,\ldots,j\}$ so that for each $v\in V(X)$, $f(N[v])\geq k$, where $N[v]$ is the closed neighbourhood of $v$. Such a function is \emph{efficient} if all of the vertex inequalities are met with equality. They give a simple necessary condition for efficient domination, namely: if $X$ is an $r$-regular graph on $n$ vertices that has an efficient $(1,k)$-dominating function, then the size of the corresponding dominating set divides $n\cdot k$.

The Hamming graph $H(q,d)$ is the graph on the vectors $\mathbb{Z}_q^d$ where two vectors are adjacent if and only if they are at Hamming distance $1$. We show that if $q$ is prime, then the previous necessary condition is sufficient for $H(q,d)$ to have an efficient $(1,k)$-dominating function. This result extends a result of Lee (Jaeun Lee. Independent perfect domination sets in Cayley graphs. \emph{J. Graph Theory}, 37(4):213–219, 2001.) on independent perfect domination in Cayley graphs. We mention difficulties that arise for $H(q,d)$ when $q$ is a prime power but not prime.
\end{abstract}


Efficient, or perfect, domination in graphs was first introduced by Biggs \cite{Biggs} as a generalization of perfect codes to graphs. Many variants of domination, and efficient domination, have been studied (see \cite{MR1605684} for a survey). Of these, $(j,k)$-dominating functions (introduced by Rubalcaba and Slater \cite{MR2412356}) give a framework that generalizes both dominating sets, and $k$-dominating sets (here $k$ refers to the number of vertices each vertex is dominated by). In this paper we consider efficient $(j,k)$-dominating functions on regular graphs.

Constructing and characterizing efficient $(j,k)$-dominating functions on a graph is challenging. For example, Rubalcaba and Slater characterize efficient $(j,k)$-domination on trees, and characterization of the efficient $(j,k)$-dominating functions on cycles follows from Lee's characterization of efficient dominating sets in Cayley graphs for Abelian groups \cite{Lee}. But there are very few classes of graphs for which we have a complete characterization.

Efficient $(j,k)$-dominating functions are related to the $(-1)$-eigenspace of a graph. When $X$ is regular, this connection is particularly strong. In Section \ref{Value} we show that a regular graph has a non-trivial efficient $(j,k)$-dominating function for some $j$ and $k$ if and only if it has $-1$ as an eigenvalue. 

 In Section \ref{DominatablePartitions}, we develop \emph{dominatable partitions} which are particularly useful for constructing efficient $(j,k)$-dominating functions. We extend covers of graphs to $m$-covers, which are a type of dominatable partition. As a proof of the usefulness of this idea, we show that Hamming graphs are $m$-covers of complete graphs.

The characterization of perfect codes by Tiet{\"a}v{\"a}inen, and Zinoviev and Leontiev \cite{Biggs2} is a classic result in coding theory. Seen as a statement about the existence of efficient dominating sets, it characterizes exactly which Hamming graphs $H(q,d)$ for $q$ a prime power admit efficient dominating sets. Efficient dominating sets are efficient $(1,1)$-dominating functions. Our main result is to partially extend this result to efficient $(1,k)$-dominating functions. When $p$ is a prime, we are able to completely characterize the efficient $(1,k)$-dominating functions on $H(p,d)$. Theorem \ref{HammingCharacterization} demonstrates that for $p$ prime, the necessary divisibility condition for an efficient $(1,k)$-dominating function is sufficient. From the proof we see that when $q$ is a prime power that is not prime, we are able to construct efficient $(1,k)$-dominating functions for many of the values of $k$ predicted by our divisibility condition. However, we fall short of a full characterization. We conjecture that the divisibility condition is always sufficient.

Finally, in Section \ref{FoldedCubes} we look at folded cubes, which are covers of the Hamming graphs $H(2,d)$ (the hypercubes). From the eigenvalue condition on the folded cubes we are able to deduce more about the form of the efficient $(j,k)$-dominating functions on $H(2,d)$.

\section{Preliminaries}\label{Preliminaries}

For a vertex $v$ in graph $X$, let $N(v)$ be the set of vertices adjacent to $v$, and $N[v]=N(v)\cup\{v\}$ be the \emph{closed neighbourhood} of $v$. Given a function $f$ defined on the vertices of $X$, and $S\subseteq V(X)$, we let $f(S)=\sum_{v\in V}f(v)$. In \cite{MR2412356}, Rubalcaba and Slater define a \emph{$(j,k)$-dominating function} on $X$ as a function $f:V(X)\rightarrow \{0,\ldots,j\}$ so that for each $v\in V(X)$, $f(N[v])\geq k$. Such a function is \emph{efficient} if all of the vertex inequalities are met with equality.

For example, in the six-cycle $C_6$, if we assign value $1$ to a pair of vertices at distance $3$, and $0$ to the other four vertices, we have an efficient $(1,1)$-dominating function. In the complete bipartite graph $K_{2,3}$, if we assign value $1$ to the vertices with degree $2$, and $2$ to the vertices with degree $3$, we have an efficient $(2,5)$-dominating function. In any complete graph, any assignment of values in $\{0,\ldots,j\}$ results in an efficient $(j,k)$-dominating function where $k$ is the total value assigned to the vertices.

\begin{figure}[h]
\centering
\begin{tikzpicture}[scale=0.5, vrtx/.style args = {#1/#2}{%
      circle, draw, fill=black, inner sep=0pt,
      minimum size=6pt, label=#1:#2}]
\node (0) [vrtx=above/1] at (0,2) {};
\node (1) [vrtx=right/0] at (1.732,1) {};
\node (2) [vrtx=right/0] at (1.732,-1) {};
\node (3) [vrtx=below/1] at (0,-2) {};
\node (4) [vrtx=left/0] at (-1.732,-1) {};
\node (5) [vrtx=left/0] at (-1.732,1)  {};
\node (11) [vrtx=left/2] at (8,1)  {};
\node (12) [vrtx=left/2] at (8,-1)  {};
\node (21) [vrtx=right/1] at (11,2)  {};
\node (22) [vrtx=right/1] at (11,0) {};
\node (23) [vrtx=right/1] at (11,-2) {};
\draw (0) edge (1);
\draw (1) edge (2);
\draw (2) edge (3);
\draw (3) edge (4);
\draw (4) edge (5);
\draw (5) edge (0);
\draw (11) edge (21);
\draw (11) edge (22);
\draw (11) edge (23);
\draw (12) edge (21);
\draw (12) edge (22);
\draw (12) edge (23);
\end{tikzpicture}
\caption{An efficient $(1,1)$-dominating function on $C_6$, and an efficient $(2,5)$-dominating function on $K_{2,3}$.}\label{SimpleExamples}
\end{figure}
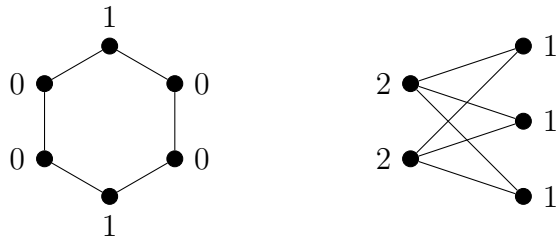

Note that from our definition, every graph has an efficient $(0,0)$-dominating function. To remove this trivial example, we require that $k\geq 1$.

Also, for any function $f:V(X)\rightarrow\{0,\ldots,j\}$, if $v$ is a vertex, and $\deg(v)=|N(v)|$, then $f(v)\leq j(1+\deg(v))$. Thus if $f$ is an efficient $(j,k)$-dominating function, we have $k\leq j(1+\delta(X))$, where $\delta(X)$ is the minimum degree of $X$. Our focus is on regular graphs of degree $r$, so this condition becomes $k\leq j(1+r)$. Moreover, we see that assigning $1$ to every vertex results in an efficient $(1,d+1)$-dominating function. When considering efficient $(1,k)$-dominating functions, we consider the all-ones function trivial.

For $r$-regular graphs, we have two elementary necessary conditions on the existence of an efficient $(1,k)$-dominating function, both of which appear in \cite{MR2412356}.
\begin{lemma}\label{divisibility}
If $X$ is $r$-regular with $n$ vertices, and $X$ has an efficient $(j,k)$-dominating function, then $r+1$ divides $n\cdot k$. Moreover, if $X$ has an efficient $(1,k)$-dominating function, then the size of the corresponding dominating set divides $n\cdot k$. \qed
\end{lemma}
The second condition pertains to efficient $(1,k)$-dominating functions.
\begin{theorem}[Theorem 10 in \cite{MR2412356}]\label{Observation}
For an $r$-regular graph $X$, and $1\leq k\leq r$, $X$ is efficiently $(1,k)$-dominatable if and only if there is a set $S\subseteq V(X)$ so the the induced subgraph $X[S]$ is $(r-1)$-regular and the induced subgraph $X[\overline{S}]$ is $(r-k)$-regular. Also $X$ is efficiently $(1,k)$-dominatable if and only if $X$ is efficiently $(1,r-k+1)$-dominatable. \qed
\end{theorem}
The restriction $1\leq k\leq r$ in this theorem only rules out the trivial dominating functions where we assign either $0$ or $1$ to all of the vertices. Assigning $0$ to all of the vertices gives $k=0$, and assigning $1$ to all of the vertices gives $k=r+1$.

\section{Linear Algebra}\label{Value}

A partition $\pi$ of $V(X)$ with cells $C_1,\ldots,C_s$ is an \emph{equitable partition} if the number of neighbours of $u\in C_i$ that lie in $C_j$ depends only on the indices $i$ and $j$. If we denote this value as $b_{ij}$, then the $s\times s$ matrix $A_{\pi}$ with $ij$-entry $b_{ij}$ is called the \emph{characteristic matrix} of the partition.

Using the vocabulary of equitable partitions, we can re-state the first part of Theorem \ref{Observation} as follows.
\begin{theorem}\label{Observation2}
For an $r$-regular graph $X$, and $1\leq k\leq r$, $X$ is efficiently $(1,k)$-dominatable if and only if $X$ has an equitable partition $\pi$ with characteristic matrix
\[
A_{\pi}=\left[\begin{array}{rr}
k-1 & r-k+1\\
k & r-k
\end{array}\right].
\]
\qed
\end{theorem}

For a more thorough treatment of equitable partitions we refer the reader to Section 9.3 of \cite{GodsilRoyle}. The following theorem connects the spectrum of $A_{\pi}$ and the spectrum of $X$.
\begin{theorem}[Theorem 9.3.3 in \cite{GodsilRoyle}]\label{evals}
If $\pi$ is an equitable partition of a graph $X$, then the characteristic polynomial of $A_\pi$ divides the characteristic polynomial of $A(X)$.

\qed
\end{theorem}
This theorem immediately implies that every eigenvalue $\theta$ of $A_{\pi}$ is an eigenvalue of $X$, and the multiplicity of $\theta$ as an eigenvalue of $A_{\pi}$ gives a lower bound on its multiplicity as an eigenvalue of $X$. The following corollary results from combining Theorems \ref{Observation2} and \ref{evals}.
\begin{corollary}\label{equitable}
For an $r$-regular graph $X$, if $X$ has an efficient $(1,k)$-dominating function for some $1\leq k\leq r$, then $-1$ is an eigenvalue of $X$. \qed
\end{corollary}

We can make the connection between efficient $(j,k)$-dominating functions and the adjacency matrix more explicit. For any assignment $f:V(X)\rightarrow\{0,\ldots,j\}$, we let the vector $\vec{f}\in\RR^n$ be the vector with $v$-component $f(v)$. Now $f$ is an efficient $(j,k)$-dominating function if and only if
\[
(A(X)+I_n)\vec{f}=k\vec{1},
\]
where $I_n$ is the $n\times n$ identity matrix. This gives us a definition of efficient $(j,k)$-domination using only linear algebra. For regular graphs this is particularly helpful.

If $X$ is $r$-regular, then every row of its adjacency matrix sums to $r$. Thus the all-ones vector, $\vec{1}$, lies in the column space of $A(X)$, and $\vec{1}\in\col(A(X)+I_n)$. This implies that the system $(A(X)+I_n)x=k\vec{1}$ is always feasible, with trivial solution $x=k/(r+1)\vec{1}$. Thus the system has non-trivial solutions if and only if $-1$ is an eigenvalue of $A(X)$. Note that only non-trivial \emph{non-negative integral} solutions qualify as efficient dominating functions. However, the necessary condition that $-1$ must be an eigenvalue is also sufficient.
\begin{theorem}\label{-1}
An $r$-regular graph $X$ has a non-trivial efficient dominating function if and only if $-1$ is an eigenvalue of $X$.
\end{theorem}
\begin{proof}
From the preceding paragraph, all that remains is to show that if $-1$ is an eigenvalue of $X$, then $X$ has some non-trivial efficient dominating function. Note that since $A(X)$ is rational-valued, and $-1$ is rational, $A(X)$ has a rational-valued $(-1)$-eigenvector. Let $x$ be an integral $(-1)$-eigenvector for $X$. Since $r$ is a simple eigenvalue for $A(X)$ with eigenspace spanned by $\vec{1}$, the $(-1)$-eigenspace is orthogonal to $\vec{1}$. Thus $x$ has both positive and negative entries.

Now
\[
(A(X)+I_n)(x+a\vec{1})=\vec{0}+a(r+1)\vec{1}=a(r+1)\vec{1}
\]
for any $a$. Choose $a$ to be the smallest positive value so that $x+a\vec{1}$ is non-negative. If $m$ is the maximum value of the components of $x+a\vec{1}$, then the result is an efficient $(m,a(r+1))$-dominating function on $X$.
\end{proof}

Theorem \ref{-1} characterizes the regular graphs that are efficiently dominatable. But it leaves something to be desired in that it tells us nothing about the values $j$ and $k$ for which efficient $(j,k)$-dominating functions exist.

Finally, we note that the work of Cardoso (and various coauthors) \cite{MR2099043,MR3546269,MR2311101,MR2599869} pushes this idea in a different direction. They look at \emph{$(k,\tau)$-partitions} of regular graphs. That is equitable partitions $(S,\overline{S})$ where $S$ induces a $k$-regular subgraph, and every vertex in $\overline{S}$ has exactly $\tau$ neighbours in $S$. Their focus is on the eigenvalues, and eigenspaces related to these partitions.

\section{Dominatable Partitions}\label{DominatablePartitions}

Here we extend the material from Section \ref{Value} towards constructing efficient $(j,k)$-dominating functions. We define a \emph{dominatable partition} of a graph $X$ to be an equitable partition $\pi$ with cells $C_1,\ldots,C_s$, and adjacency constants $b_{il}$ satisfying the following condition. For each $1\leq l\leq s$ there is a value $a_l$ so that $b_{il}=a_l$ for all $i\neq l$, and $b_{ll}=a_l-1$. Combinatorially this means that every vertex outside of $C_l$ is adjacent to exactly $a_l$ vertices in $C_l$, and every vertex in $C_l$ is adjacent to exactly $a_l-1$ vertices in $C_l$. Note that the degree of any vertex of $X$ is $r=(\sum a_l)-1$, and $X$ is $r$-regular.

Given such a partition, and values $\alpha_l$ for $0\leq l\leq s$, we define $f:V(X)\rightarrow\RR$ as $f(v)=\alpha_l$ for all $v\in C_l$. Now for any vertex $v\in V(X)$, we have
\[
f(N[v])=\sum_{l=1}^s\alpha_la_l=k.
\]
Note that $k$ is constant. If we take $\alpha_l\in\{0,\ldots,j\}$ for all $l$, we have an efficient $(j,k)$-dominating function.

If we look at the characteristic matrix $A_{\pi}$ of a dominatable partition, we see that $A_{\pi}+I_s$ has constant columns. Thus we can express this matrix as an outer product
\[
A_{\pi}+I_{s}=\vec{1}\cdot[a_1\ldots a_s].
\]
Therefore, $\vec{1}$ is an eigenvector with corresponding eigenvalue $\sum a_l=r+1$, and the remaining $s-1$ eigenvalues are all $0$. Thus $A_{\pi}$ has eigenvalue $r$ with multiplicity $1$, and eigenvalue $-1$ with multiplicity $s-1$. We immediately have the following analogue of Corollary \ref{equitable}.
\begin{theorem}\label{dominatable}
If $X$ is an $r$-regular graph with dominatable partition $\pi$ with $s$ parts, then $-1$ is an eigenvalue of $X$ with multiplicity at least $s-1$.\qed
\end{theorem}

Note that Theorem \ref{dominatable} does not bring us much closer to a characterization of the pairs $(j,k)$ for which an efficient dominating function exists. But if we were instead to ask for efficient $(j,k)$-dominating functions that take $s$ distinct values, it gives us a condition that rules out such dominating functions arising from dominatable partitions. It also is not a characterization of graphs with dominatable partitions.

We can also characterize dominatable partitions by their eigenvalues. Suppose $\pi$ is an equitable partition of an $r$-regular graph $X$ with $s$ parts whose characteristic matrix $A_{\pi}$ has eigenvalue $-1$ with multiplicity $s-1$. Since $X$ is $r$-regular, $A_{\pi}$ has eigenvalue $r$ with multiplicity $1$. So $A_{\pi}+I_s$ has eigenvalue $r+1$ with multiplicity $1$ and eigenvalue $0$ with multiplicity $s-1$. Since the row sums of $A_{\pi}+I_s$ are constant, this immediately implies that this matrix is an outer product of the form
\[
A_{\pi}+I_s=\vec{1}[a_1\ldots a_r].
\]
Thus $\pi$ is a dominatable partition.

As an example, consider the complete graph $K_n$. The partition of $V(K_n)$ into singletons is equitable with characteristic matrix $A(K_n)$. The eigenvalues of $K_n$ are $n$ with multiplicity $1$ and $-1$ with multiplicity $n-1$, so this partition is dominatable.

In Section \ref{HammingGraphs} we will see that the $5$-cube is an example of a graph that has $-1$ as an eigenvalue with multiplicity larger than one, but has no dominatable partition of size more than $2$. This means that the equitable partitions of the $5$-cube with $-1$ as an eigenvalue are either partitions with $2$ parts, or partitions that have more than $2$ eigenvalues. We will also see that the $3$-cube has many dominatable partitions. This is because $Q_3$ is a cover of $K_4$.

\section{Covers}\label{Covers}

A graph $X$ is a \emph{cover} of a graph $Y$ if there is a partition $\{C_v\,:\,v\in V(Y)\}$ of $V(X)$ so that: each $X[C_v]$ is independent; and, $X[C_v\cup C_u]$ has no edges when $u$ and $v$ are non-adjacent in $Y$, and is a perfect matching when $u$ and $v$ are adjacent. From the definition we see that each $C_v$ must have the same size. If each $|C_v|=a$ then we say $X$ is an \emph{$a$-fold cover} of $Y$. The sets $C_v$ are called the \emph{fibres} of the cover.

Note that if $X$ is a cover of $K_m$, then the partition given by the fibres is a dominatable partition with characteristic matrix $J_m-I_m$. If $X$ is $r$-regular, then $m=r+1$, and $r+1$ must divide $n$ the number of vertices of $X$. If we assign $1$ to one cell, and $0$ to the rest, we obtain an efficient $(1,1)$-dominating function. In \cite{Lee}, Lee proves that if $X$ is a Cayley graph for an Abelian group, then it has an efficient $(1,1)$-dominating function if and only if it is a cover of a complete graph (the language used by Lee is an ``independent perfect domination set''). We draw the statement of the following theorem from the results in \cite{Lee}.
\begin{theorem}\label{LeesTheorem}
If $X=X(G,C)$ is a Cayley graph for Abelian group $G$, and $S$ is the set of vertices assigned $1$ by an efficient $(1,1)$-dominating function on $X$, then:
\begin{enumerate}
\item for each $c\in C$, the function $f_c$ that assigns value $1$ to the vertices in $cS$, and $0$ otherwise, is an efficient $(1,1)$-dominating function on $X$;
\item the sets $\{S\}\cup\{cS\,:\,c\in C\}$ partition $G$;
\item the map $\rho: G\rightarrow V(K_{|C|+1})$ where each set from (2) is mapped to a distinct vertex of $K_{|C|+1}$ is a covering map.
\end{enumerate}\qed
\end{theorem}
Note that in general we only require the edges between fibres of a cover to be a perfect matching. But in a Cayley graph, each edge is associated with a generator in $C$, and in the cover described in Theorem \ref{LeesTheorem} the matching connecting two fibres consists of edges generated by a single element of $C$. We will make use of this additional property when considering Hamming graphs. In our application we will additionally start with a set $S$ that is a subgroup of $G$, whence the fibres of the cover will be the cosets of $S$. Finally note that Lee's result proves that in a Cayley graph for an Abelian group, every efficient $(1,1)$-dominating function arises from a dominatable partition.

We can also consider generalizing the concept of a cover to a partition that leads to the construction of efficient dominating functions. We define $X$ to be a \emph{$k$-cover} of a graph $Y$ if there is a partition $\{C_v\,:\, v\in V(Y)\}$ of $V(X)$ so that: each $X[C_v]$ is a $(k-1)$-regular graph; and, $X[C_u\cup C_v]$ has no edges when $u$ and $v$ are non-adjacent in $Y$, and is a $k$-regular bipartite graph when $u$ and $v$ are adjacent. Note that the partition given by a $k$-cover of $K_m$ is dominatable.

We noted that if $X$ is a cover of $K_m$, then that cover leads to the construction of an efficient $(1,1)$-dominating function on $X$. Likewise, if $X$ is a $k$-cover of $K_m$, then by assigning one of the fibres $1$, and the rest $0$, we have an efficient $(1,k)$-dominating function on $X$. As a trivial example of such a partition, we can take an arbitrary partition of the vertices of $K_{2n}$ into sets of size $2$. This gives a $2$-cover of $K_n$ by $K_{2n}$. In Section \ref{HammingGraphs} we will see non-trivial examples.

\section{Hamming Graphs}\label{HammingGraphs}

Hamming graphs are a classical setting for domination problems. Given a set $\Omega$ of $q$ symbols, and a length $d\geq 1$, the \emph{Hamming graph} $H(q,d)$ is the graph on $\Omega^d$ where two tuples are adjacent if and only if they differ in exactly one component. Alternatively, they are adjacent if and only if they are at Hamming distance one. Note that $H(q,d)$ is a graph on $q^d$ vertices, and every vertex has exactly $(q-1)d$ neighbours. If $q=2$, then $H(2,d)=Q_d$ is the $d$-dimensional hypercube.

Our focus will mostly be on Hamming graphs $H(q,d)$ where $q$ is a prime power. In this case we can take the symbols $\Omega$ to be the elements of the finite field of order $q$. So we can describe our Hamming graph as the Cayley graph $H(q,d)=X(GF(q)^d,C_d)$ where $\mathcal{C}_d=\{\alpha\vec{e}_i\,:\,1\leq i\leq d,\ \alpha\neq 0\}$ (here $\vec{e}_i$ is the $i$th standard basis vector).

Distance in $H(q,d)$ corresponds to the Hamming distance between tuples in $\Omega^d$, and thus are of interest to coding theorists. A \emph{perfect $e$-code} in a graph $X$ is a set $S\subseteq V(X)$ so that that every vertex of $X$ is at distance at most $e$ from exactly $1$ member of $S$. This concept was introduced by Biggs in \cite{Biggs}. Efficient $(j,k)$-dominating functions can be viewed as a generalization of perfect $1$-codes, as a perfect $1$-code is exactly an efficient $(1,1)$-dominating function. For Hamming graphs, perfect $e$-codes have a complete characterization due independently to Tiet{\"a}v{\"a}inen, and to Zinoviev and Leontiev \cite{Biggs2}. In particular, the only Hamming graphs that admit perfect $1$-codes are $H(q,(q^a-1)/(q-1))$ where $q$ is a prime power. The linear perfect $1$-codes in these graphs are the Hamming codes. Our main result in this section is the extension of this theorem from efficient $(1,1)$-dominating functions to efficient $(1,k)$-domination functions.
\begin{theorem}\label{HammingCharacterization}
Let $q$ be a prime power. If $(q-1)d+1=q^a m$ where $m$ is not divisible by $q$, then $H(q,d)$ has an efficient $(1,k)$-dominating function for all $1\leq k\leq (q-1)d+1$ a multiple of $m$. Moreover, if $q$ is prime, then $H(q,d)$ has an efficient $(1,k)$-dominating function if and only if $1\leq k\leq (q-1)d+1$ is a multiple of $m$.
\end{theorem}
Note that Theorem \ref{HammingCharacterization} says that the necessary divisibility condition given in Lemma \ref{divisibility} is sufficient for the existence of efficient $(1,k)$-dominating functions for primes $p$. Before giving the full proof, it is helpful to consider two extremes.

If $(q-1)d+1$ is not divisible by $q$, then our necessary condition states that the only values $0\leq k\leq (q-1)d+1$ for which $H(q,d)$ has an efficient $(1,k)$-dominating function are $k=0$ and $k=(q-1)d+1$. The trivial efficient dominating functions, the zero function and the all-ones function, give us the efficient $(1,k)$-dominating functions in this case. So the necessary condition is trivially sufficient.

If $(q-1)d+1=q^a$ for some $a$, then the necessary condition allows for the existence of an efficient $(1,1)$-dominating function. There is a Hamming code $C$ over $GF(q)$ with length $(q^a-1)/(q-1)$, dimension $(q^a-1)/(q-1)-a$ and distance $3$. Moreover, this code is perfect. That is, that there is a subspace $C$ of $GF(q)^l$ with $l=(q^a-1)/(q-1)$ so that the function that assigns $1$ to the elements of $C$ and $0$ to the elements outside of $C$ gives an efficient $(1,1)$-dominating function on $H(q,d)$. Moreover, from Theorem \ref{LeesTheorem}, the partition given by the cosets of $C$ are the fibres of a cover of $K_{d+1}$ by $H(q,d)$, and form a dominatable partition. Thus the existence of Hamming codes proves that our necessary condition is sufficient when $(q-1)d+1=q^a$ for some $a$.

The Hamming codes in the preceding paragraph will be an important base case in our proof, so we should address the case when $a=1$. If $a=1$, then $d=1$ and we easily find our efficient dominating functions as $H(q,1)=K_q$. We want an analogue for Hamming codes, so by convention we take $\{\vec{0}\}$ to be our code when $a=1$. This is the singleton code in $H(q,1)$, it is also a subspace with basis $\emptyset$. The most important property of this code is that it is perfect, as is the case for the Hamming codes.

With these preliminaries out of the way, we can now prove our full characterization.
\begin{proof}[Proof of Theorem \ref{HammingCharacterization}]
From the preceding comments it suffices to prove the claim when $(q-1)d+1=q^am$ for $a\geq 1$ and $m>1$. We do this by constructing a dominatable partition of $H(q,d)$ that is an $m$-cover of $K_{q^a}$. Throughout the proof we take $l=(q^a-1)/(q-1)$.

Let $\vec{f}_i$ be the $i$th standard basis vector in $GF(q)^l$, and $\vec{e}_j$ be the $j$th standard basis vector in $GF(q)^d$. Partition the vectors $\vec{e}_j$ into sets $S_i$ indexed by $0\leq i\leq l$ arbitrarily so that $|S_0|=(m-1)/(q-1)$, and $|S_i|=m$ for $i>0$. This is possible as $(q-1)d=q^am+1$, which can be rearranged to
\[
d=\frac{q^am-1}{q-1}=\left(\frac{q^a-1}{q-1}\right)m+\frac{m-1}{q-1}=lm+\frac{m-1}{q-1}.
\]
Thus $(m-1)/(q-1)$ is an integer, and we can for a partition of our $d$ standard basis vectors $\vec{e}_j$ into the sets $S_i$ as described.

We use this partition to define a linear transformation $\varphi:GF(q)^d\rightarrow GF(q)^l$. For each $\vec{e}_j$ we map
\[
\varphi(\vec{e}_j)=\begin{cases}
\vec{f}_i & \text{if $\vec{e}_j \in S_i$ and $1\leq i\leq l$,}\\
\vec{0} & \text{if $\vec{e}_j \in S_0$}.
\end{cases}
\]
Then we extend $\varphi$ linearly to all of $GF(q)^d$.

Using this map, we define the promised dominatable partition of $H(q,d)$. Let $C$ be a perfect code in $GF(q)^l$ (either a Hamming code if $a\geq 2$, or the singleton code if $a=1$), and let $B_C$ be a basis for $C$. Let $T\subseteq GF(q)^d$ be the pre-image of $C$ under $\varphi$. Then $T$ is a subspace, and to find $\dim(T)$, we construct a basis.

First, we find a basis $B$ of $\ker(\varphi)$. Note that if $v\in GF(q)^d$ is mapped to $\vec{0}$, then
\[
\vec{0}=\varphi(v)=\varphi\left(\sum_{j=1}^dv_j\vec{e}_j\right)=\sum_{i=1}^dv_j\varphi(\vec{e}_j)=\sum_{i=1}^l\left(\sum_{\vec{e}_j\in S_i}v_j\right)\vec{f}_i
\]
(as $\varphi(\Span(S_0))=\{\vec{0}\}$). We see that $\varphi(v)=\vec{0}$ if and only if
\begin{equation}\label{subspace}
\sum_{\vec{e}_j\in S_i}v_j=0
\end{equation}
for each $1\leq i\leq l$.

For each $1\leq i\leq l$, let $B_i$ be a basis for the subspace of $\Span(S_i)$ given by the vectors that satisfy Equation \eqref{subspace}. Since $\dim(\Span(S_i))=m$, $|B_i|=m-1$. Finally, we take $B_0$ to be a basis for $\Span(S_0)$, and note $|B_0|=(m-1)/(q-1)$. Taking $B=\cup_{i=0}^l B_i$ gives a basis for $\ker(\varphi)$, and
\[
|B|=l(m-1)+\frac{m-1}{q-1}=\frac{q^a(m-1)}{q-1}.
\]

For each $b\in B_C$ (our basis for $C$) take an arbitrary element $b_{\varphi}$ in the preimage of $b$ under $\varphi$. Now
\[
B'=B\cup\{b_{\varphi}\,:\,b\in B_C\}
\]
is a basis for $T$. Note that
\[
|B'|=|B|+|B_C|=\frac{q^a(m-1)}{q-1}+\frac{q^a-1}{q-1}-a=d-a.
\]
Thus $\dim(T)=d-a$, and $T$ has $q^a$ cosets in $GF(q)^d$. Moreover, $\varphi$ gives us a bijection between the cosets of $T$ in $GF(q)^d$ and the cosets of $C$ in $GF(q)^l$.

Consider the partition of $H(q,d)$ given by the cosets of $T$. Since the linear transformation $\rho_v(x)=x+v$ on $GF(q)^d$ is an automorphism of $H(q,d)$, the subgraphs induced by each coset of $T$ are isomorphic. Let $T_u$ and $T_v$ be any two distinct cosets of $T$. Under $\varphi$, these cosets correspond to two distinct cosets $C_{u'}$ and $C_{v'}$ of $C$. Since the cosets of $C$ give a cover of $K_{q^a}$, the edges between $C_{u'}$ and $C_{v'}$ in $H(q,l)$ form a perfect matching. Morever, this perfect matching corresponds to a single generator $g\vec{f}_i$. Thus in $H(q,d)$, the edges between $T_a$ and $T_b$ all correspond to the generators in $gS_{i}$. Therefore the edges between $T_a$ and $T_b$ form an $m$-regular bipartite graph. Finally, we note that for $v\in T$, the number of edges between $v$ and elements of the other cosets of $T$ is
\[
m(q^a-1)=(q-1)d-(m-1)
\]
so every coset induces an $(m-1)$-regular subgraph of $H(q,d)$. Thus the partition forms an $m$-cover of $K_{q^a}$, as required.

\end{proof}

From Theorem \ref{HammingCharacterization} we make the following conjecture.
\begin{conjecture}\label{HammingConjecture}
Let $q=p^b$ for $p$ prime. If $(d-1)q+1=p^am$ where $m$ is not divisible by $p$, then $H(q,d)$ has an efficient $(1,k)$-dominating function if and only if $1\leq k\leq (q-1)d+1$ is a multiple of $m$.
\end{conjecture}
That is, we conjecture that the necessary divisibility condition from Lemma \ref{divisibility} is sufficient for $H(q,d)$ with $q$ a prime power.

The first open case is $q=4$. For $d-1$ not divisible by $q$, $H(q,d)$ does not have $-1$ as an eigenvalue, so there are no non-trivial efficient dominating functions. Checking the first few values of $d\equiv_4 1$ we find $m$-covers following the same construction as in Theorem \ref{HammingCharacterization}. The results are summarized in the following table.
\[
\begin{tabular}{c|c|c}
$d$ & $3d+1$ & partitions\\
\hline
\hline
$5$ & $4^2\cdot 1$ & cover of $K_{16}$\\ 
$9$ & $4^1\cdot 7$ & $7$-cover of $K_4$\\
$13$ & $4^1\cdot 10$ & $10$-cover of $K_4$
\end{tabular}
\]
For $H(4,5)$ and $H(4,9)$, the $m$-cover technique from Theorem \ref{HammingCharacterization} gives all of the predicted efficient $(1,k)$-dominating functions. However for $H(4,13)$, the divisibility condition is implies that $k$ must be divisible by $5$, while the $10$-cover of $K_4$ only gives efficient $(1,k)$-dominating functions for $k\in\{10,20,30\}$. We still do not know whether there are efficient $(1,k)$-dominating functions for $k\in\{5,15,25,35\}$.

\section{Folded Cubes}\label{FoldedCubes}

The $d$-dimensional hypercube is the graph $Q_d=H(2,d)$. From Theorem \ref{HammingCharacterization} we have a full characterization of the efficient $(1,k)$-dominating functions for hypercubes. An interesting related family of graphs is the folded cubes. The \emph{folded cube of order $d$} is $F_d=X(\ZZ_2^{d-1},C)$ where $C=\{\vec{e}_1,\ldots,\vec{e}_{d-1},\vec{1}\}$. From this definition, $F_d$ can be constructed from a hypercube in two ways. First, $F_d$ is the graph obtained from $Q_{d-1}$ by adding a perfect matching joining each vertex to its antipode (i.e., adding the edges joining $v$ to $v+\vec{1}$). Alternatively we can construct $F_d$ by ``folding'' $Q_d$. That is, $F_d$ is the graph on the pairs of antipodal vertices where two pairs $\{v,v+\vec{1}\}$ and $\{u,u+\vec{1}\}$ are adjacent if and only if the Hamming distance from $v$ to either $u$ or $u+\vec{1}$ is $1$. This construction shows that $Q_d$ is a $2$-fold cover of $F_d$.

To connect efficient dominating functions on $F_d$ to efficient dominating functions on $Q_d$, we make the following easy observation. If $X$ is an $m$-fold cover of $Y$, and $f:V(X)\rightarrow\RR$ is constant on the fibres of the cover, then we can define $\check{f}:V(Y)\rightarrow\RR$ as $\check{f}(v)=f(u)$ where $u$ is any vertex in $C_v$. Likewise, if $f:V(Y)\rightarrow\RR$ is any function, then we can define $\hat{f}:V(X)\rightarrow\RR$ as $\hat{f}(u)=f(v)$ where $u\in C_v$. From the definitions, we have the following observation.
\begin{lemma}\label{CoverFcn}
If $X$ is an $m$-fold cover of $Y$, then any efficient $(j,k)$-dominating function $f$ on $X$ that is constant on the fibres of the cover gives an efficient $(j,k)$-dominating function $\check{f}$ on $Y$. Similarly, any efficient $(j,k)$-dominating function on $Y$ gives an efficient $(j,k)$-dominating function $\hat{f}$ on $X$. \qed
\end{lemma}
Finally, we note that from Theorem \ref{-1}, $-1$ must be an eigenvalue of $F_d$ if it has a non-trivial efficient dominating function. Since $F_d$ has $-1$ as an eigenvalue only when $d+1$ is divisible by $4$, we have the following corollary to Theorem \ref{HammingCharacterization}.
\begin{corollary}\label{QdFd}
If $d+1$ is not divisible by $4$, then no non-trivial efficient dominating function on $Q_d$ is constant on all pairs of antipodal vertices. \qed
\end{corollary}

\begin{center}
\bf{Acknowledgements}
\end{center}

The author thanks Gary MacGillivray for several helpful discussions.

\bibliographystyle{plain}
\bibliography{EfficientDomination}

\begin{thebibliography}{10}

\bibitem{Biggs}
Norman Biggs.
\newblock Perfect codes in graphs.
\newblock {\em J. Combin. Theory Ser. B}, 15:289--296, 1973.

\bibitem{Biggs2}
Norman Biggs.
\newblock Perfect codes and distance-transitive graphs.
\newblock In {\em Combinatorics ({P}roc. {B}ritish {C}ombinatorial {C}onf.,
  {U}niv. {C}oll. {W}ales, {A}berystwyth, 1973)}, number~13 in London Math.
  Soc. Lecture Note Ser., pages 1--8. Cambridge Univ. Press, London, 1974.

\bibitem{MR2099043}
D.~M. Cardoso and P.~Rama.
\newblock Equitable bipartitions of graphs and related results.
\newblock {\em J. Math. Sci. (N. Y.)}, 120(1):869--880, 2004.

\bibitem{MR3546269}
Domingos~M. Cardoso, Vadim~V. Lozin, Carlos~J. Luz, and Maria~F. Pacheco.
\newblock Efficient domination through eigenvalues.
\newblock {\em Discrete Appl. Math.}, 214:54--62, 2016.

\bibitem{MR2311101}
Domingos~M. Cardoso and Paula Rama.
\newblock Spectral results on regular graphs with {$(k,\tau)$}-regular sets.
\newblock {\em Discrete Math.}, 307(11-12):1306--1316, 2007.

\bibitem{MR2599869}
Domingos~M. Cardoso, Irene Sciriha, and Cheryl Zerafa.
\newblock Main eigenvalues and {$(\kappa,\tau)$}-regular sets.
\newblock {\em Linear Algebra Appl.}, 432(9):2399--2408, 2010.

\bibitem{GodsilRoyle}
Chris Godsil and Gordon Royle.
\newblock {\em Algebraic {G}raph {T}heory}.
\newblock Number 207 in Graduate Texts in Mathematics. Springer-Verlag, New
  York, 2001.

\bibitem{MR1605684}
Teresa~W. Haynes, Stephen~T. Hedetniemi, and Peter~J. Slater.
\newblock {\em Fundamentals of domination in graphs}, volume 208 of {\em
  Monographs and Textbooks in Pure and Applied Mathematics}.
\newblock Marcel Dekker, Inc., New York, 1998.

\bibitem{Lee}
Jaeun Lee.
\newblock Independent perfect domination sets in {C}ayley graphs.
\newblock {\em J. Graph Theory}, 37(4):213--219, 2001.

\bibitem{MR2412356}
Robert~R. Rubalcaba and Peter~J. Slater.
\newblock Efficient {$(j,k)$}-domination.
\newblock {\em Discuss. Math. Graph Theory}, 27(3):409--423, 2007.

\end{thebibliography}

\end{document}